\documentstyle[twoside,12pt]{article}
\def\be{\begin{equation}}\def\ee{\end{equation}}
\def\benon{\begin{displaymath}}\def\eenon{\end{displaymath}}
\def\bea{\begin{eqnarray}}\def\eea{\end{eqnarray}}
\def\bma{\left(\begin{array}}\def\ema{\end{array}\right)}
\newenvironment{smallarray}{\arraycolsep1mm
                            \begin{array}}{\end{array}}
\def\bsma{\left(\begin{smallarray}}\def\esma{\end{smallarray}\right)}
\makeatletter
\newenvironment{xequation}[1]% EQUATION WITH EXPLICIT GIVEN LABEL
{\def\@eqnnum{{\bf (#1)}}\be}{\ee\addtocounter{equation}{-1}}
\newenvironment{xequationarray}[1]% EQUATIONARRAY WITH EXPLICIT GIVEN LABEL
{\def\@eqnnum{{\bf (#1)}}\be\begin{array}{rcl}}%
{\end{array}\ee\addtocounter{equation}{-1}}
\makeatother
\def\bex#1{\begin{xequation}{#1}}\def\eex{\end{xequation}}
\def\beax#1{\begin{xequationarray}{#1}}\def\eeax{\end{xequationarray}}
\newlength{\partcolwidth}\newlength{\restcolwidth}
\def\partcolumn#1#2{\setlength{\partcolwidth}{#1}
                    \raisebox{\baselineskip}{\parbox[t]{\partcolwidth}{#2}}}
\def\restcolumn#1{\setlength{\restcolwidth}{\textwidth}
                  \addtolength{\restcolwidth}{-\partcolwidth}
                  \raisebox{\baselineskip}{\parbox[t]{\restcolwidth}{#1}}}

\def\c{\cite}
\def\ii{\vphantom{)}}
\def\q#1{\ \ \hbox{#1}\ \ }
\def\qq#1{\qquad\hbox{#1}\qquad}
\def\o#1{\makebox[1em]{$#1$}}
\def\m{\makebox[0mm][r]{$-$}}
\def\nn{\nonumber}

\def\Rhat{\hat{R}}
\def\t#1{\widetilde{#1}}
\def\bino#1#2{{#1\choose #2}}
\def\dfrac#1#2{ {\displaystyle\frac{#1}{#2}} } % force big fraction
\def\tfrac#1#2{ {\textstyle\frac{#1}{#2}} }    % force small fraction

\def\inbar{\,\vrule height1.5ex width.4pt depth0pt}
\def\QQ{\relax\,\hbox{$\inbar\kern-.3em{\rm Q}$}}
\def\CC{\relax\hbox{$\inbar\kern-.3em{\rm C}$}}
\def\RR{\relax{\rm I\kern-.18em R}}
\def\ZZ{\relax{\rm Z\kern-.4em Z}}
\def\DD{{I\!\!D}}
\def\oo{\mbox{$1\!\! 1$}}
\def\a{\alpha} \def\b{\beta} \def\g{\gamma} \def\dd{\delta} \def\l{\lambda}
\def\e{\varepsilon} \def\si{\sigma} \def\k{\kappa} \def\z{\zeta}
\def\ot{\!\otimes\!}
\def\cd{\!\cdot\!}
\def\dpo{{\displaystyle 1}}    % for use in exponents
\def\tr{\hbox{\rm tr}}
\def\S{\scriptstyle}

\begin{document}

\thispagestyle{empty}
\ii\vspace{-2cm}
\rightline{LMU-TPW 94/23}
\rightline{MPI-PhT/94-93}
\rightline{December 1994}
\vspace{2ex}
\begin{center}
{\bf \LARGE Classification of the GL(3) Quantum Matrix Groups}\\[6ex]
Holger\, Ewen\footnote{email: hoe@mppmu.mpg.de}\\[2ex]
Sektion Physik, Universit\"at M\"unchen\\
Theresienstr. 37, D-80333 M\"unchen, Germany\\[4ex]
Oleg\, Ogievetsky\footnote{on leave of absence from N. Lebedev Physical
Institute, Theoretical Department, Leninsky prospect~53, 117924 Moscow,
Russia}\\[2ex]
Max--Planck--Institut f\"ur Physik\\
F\"ohringer Ring 6, D-80805 M\"unchen, Germany \\[6ex]
\end{center}

{\bf Abstract.}
We define quantum matrix groups GL(3) by their coaction on appropriate
quantum planes and the requirement that the Poincar\'e series coincides with
the classical one. It is shown that this implies the existence of a Yang-Baxter
operator. Exploiting stronger equations arising at degree four of the
algebra, we classify all quantum matrix groups GL(3).
We find 26 classes of solutions, two of which do not admit a normal ordering.
The corresponding R-matrices are given.

\section*{Introduction}
In addition to the standard deformation of GL(3), in the past a few
nonstandard R-matrices were discovered which define quantum groups of GL(3)
type. These are for example the multiparameter
deformation \c{Sud,FeiGuo} and the Jordanian deformation \c{DMMZ}.
Both turn out to be related to the standard deformation.
The former is connected to it by the twisting procedure \c{Res},
but requires the knowledge of the appropriate twisting operator.
The latter is a limiting case of the standard deformation \c{Ogi2}.
Lacking a more constructive procedure, attempts have been made to find
the general solution to the classification problem by `equation
crunching', using the Yang-Baxter equation as a starting point.
This has been done successfully in 2 dimensions \c{Hie}, but is feasible in 3
dimensions only with the special ansatz of upper triangularity for the
R-matrix \c{Hie2}. These methods give R-matrices of all possible
projector decompositions, whereas for GL(3) the R-matrix has to be a sum
of two complementary projectors.
They give rise to the definition of even and
odd quantum planes whose coordinates have commutation relations determined
by the two projectors. The rank of the projectors
corresponds to the dimension of these coordinate algebras at degree two.

In the case of GL(3) it is natural to define a quantum matrix group
not by a solution of the Yang-Baxter equation but by its (co-) action on the
appropriate quantum planes, together with the requirement that the
Poincar\'e series of the algebra generated by the matrix elements
coincides with that of the classical case.
It will be shown that these demands imply the existence of an R-matrix for the
algebra, which satisfies the Yang-Baxter equation.

A three dimensional quantum plane is given as the polynomial algebra
generated by three elements which obey quadratic relations, such that the
Poincar\'e series of the algebra, given by the dimensions of the subspaces
of specific degree, coincides with the classical case of 3 commuting
generators. The consistency of a left action of a quantum matrix on a
quantum plane, together with a right action on a quantum coplane, is well
known to impose a compatibility condition on plane and coplane and to fix the
relations of the group completely. The above mentioned odd plane is defined
by the relations dual to those of the coplane.
Imposing the condition of classical Poincar\'e series on the group algebra
results in additional compatibility conditions for plane and coplane.
At degree three we get back the Yang-Baxter equation, but we find more
conditions at degree four.
One implication of these necessary conditions is that the characteristic
matrices of
plane and coplane, which we will define in section one, have to commute.
This generalizes to higher dimensions.
It is known from counterexamples that the Yang-Baxter equation is not
sufficient to fix the dimension of the algebra at all degrees to the
classical values.
The stronger equations we get allow us, in three dimensions, to solve
the classification problem case by case.
In two dimensions it is even possible to solve the equations algebraically
and to get straightforwardly the complete answer, as we have shown in
\c{EOW,EOW2}.
By following this algebraic approach we forget about quantum groups as
deformations
of classical groups. Instead of demanding continuity of the algebra
dimensions during deformation we fix them by hand. Having found a
solution, we may ask if a classical limit is possible, anyhow we will not
do this here.

In most cases the classical Poincar\'e series is accompanied by the
existence of a complete and unique ordering prescription in the sense of
the diamond lemma. We find two related cases of quantum groups for which
this is
not true, however we are not able to prove for these the correctness of the
full Poincar\'e series. We will come back to this question in section five.
If $q$, the negative quotient of the eigenvalues of the R-matrix,
is not a root of unity, then all of our equations are consequences
of the Yang-Baxter equation, as follows from Gurevich \c{Gur}.
In this restricted case the Poincar\'e series of all algebras mentioned
above is in fact determined by its first three elements.
Many of our solutions however correspond to $q$ being a root of unity.

Sections 1 and 2 derive necessary conditions for the correct
Poincar\'e series of the quantum plane and group, respectively.
Section 1 also gives the classification of three-dimensional quantum
planes by their characteristic matrix.
Another significant quantity commuting with the
characteristic matrices of the planes is the matrix which gives the
commutation law of
the quantum determinant with the group generators.
Section 3 introduces a transformation on the space of solutions
which preserves the Poincar\'e series and which we use to
reduce the number of cases for the determinant commutation
law to a few.
Starting from a specific case, we find in section 4 the explicit
solutions to the quantum group conditions.
At the end of that section we summarize their
automorphism groups, leading back to the complete solution space of the
problem.
Section 5 returns to the question of the correct Poincar\'e series
and proves correctness to all orders by giving an appropriate normal
ordering for the generators.
Appendix A lists the explicit R-matrices for our solutions.
A kind of summary of our results is given in table 3 at the beginning of
section~4.

\section{Quantum planes in three dimensions}

Quantum planes of dimension three possessing the standard Poincar\'e series
have been classified in \c{ArtSch}. However, several cases which
we would like to distinguish have not been treated there as separate cases.

Let the quadratic relations on the coordinates of the quantum plane
be given by
\be E^\a_{ij}x^ix^j=0\qq, i,j=1,2,3\label{3.1}\ee
(where summation over paired indices is understood),
which we will write in  abbreviated form as $E^\a x\cd x=0$.
The condition of standard Poincar\'e series of the algebra generated by $x^i$
requires three independent relations (\ref{3.1}),
i.e.~$\a=1,2,3$ and $E^\a$ linearly independent.
In the following we take Greek indices to be running over 1, 2, 3.
The index $\a$ can be transformed with an arbitrary
invertible matrix without changing the ideal generated.

Since we assume that all relations of higher degree are generated by the
quadratic relations (\ref{3.1}), the relations in degree three are
\be (E^\a\ot\oo^i)\,x\cd x\cd x = 0 \qq{and}
(\oo^i\ot E^\a)\,x\cd x\cd x = 0\label{3.2}\ee
(to be read as $E^\a_{mn}x^mx^nx^i=0$ and $E^\a_{mn}x^ix^mx^n=0$,
$(\oo^i)_j=\dd^i_j$).
These 18 relations cannot be independent.
Classically, $\bino{d+n-1}n$ of the $n$--th order products of $d$ commuting
variables are independent, resulting in 10 cubics for $d=3$.
With $3^3=27$ quantities $x\cd x\cd x$ we therefore need 17 relations and
there must be exactly one intersection in (\ref{3.2}).
Since the $E^\a$ are linearly independent, it has to be of the form
\be e_{i\a}\,(\oo^i\ot E^\a) = f_{\a i}\,(E^\a\ot\oo^i) =: E \ ,\label{3.3}\ee
$e$ and $f$ nonvanishing.
The (up to rescaling) unique solution of (\ref{3.3}) defines the tensor
$E=(E_{ijk})$, which is the deformed $\e$ tensor of the three-dimensional
quantum plane.

In case that $e$ and $f$ give nondegenerate mappings between the three
dimensional vector spaces indexed with Latin and Greek indices respectively,
eq.~(\ref{3.3}) gives a cyclicity property for the tensor $E$.
Comparing $E_{ijk}=E^\a_{ij} f_{\a k}$ and $E_{lij}=e_{l\a} E^\a_{ij}$
we get
\be E_{ijk}=Q^l_k E_{lij} \qq{,\qquad with} Q^i_j = f_{\a j} (e^{-1})^{\a i}\ .
\label{3.8}\ee
We will assume in the following that both $e$ and $f$ are nondegenerate and
will characterize quantum planes as solutions of the cyclicity equation
(\ref{3.8}).
The fact that degenerate $e$ or $f$ can also lead to  standard
Poincar\'e series  is demonstrated by the example $E^1=e^3e^2$,
$E^2=e^3e^1$, $E^3=e^2e^1$.
A matrix $Q$ may also be defined in the case where
$e$ and $f$ have equal rank and the kernels of the two mappings
$x^\a\mapsto x^\a e_{i\a}$ and $x^\a\mapsto x^\a f_{\a i}$
coincide. It can be chosen to be invertible as well.

Choosing a basis of the algebra where the invertible matrix $Q$ has Jordan
normal form, we get three cases:
\be \bma{ccc} \o{\a_1}&&\\&\o{\a_2}&\\&&\o{\a_3}\ema \qq,
    \bma{ccc} \a&1&\\&\a&\\&&\b \ema \qq,
    \bma{ccc} \a&1&\\&\a&1\\&&\a \ema\ .\label{3.9}\ee
Starting with diagonalizable $Q$, the cyclicity equation becomes
$E_{ijk}=\a_k E_{kij}$.
This gives a set of 11 components $E_{ijj}$, $E_{123}$ and $E_{132}$
determining all remaining ones and independently subject to conditions
\be E_{iii}=\a_iE_{iii} \qq, E_{ijk}=\a_i\a_j\a_kE_{ijk}\ .\ee
Hence either distinct products of the eigenvalues of $Q$ have to be equal
to one or the corresponding components of the $E$ tensor vanish.
The different possibilities for the eigenvalues of $Q$ give rise to the
cases 1 to 11 in the table below.
We excluded cases with dependent $E^\a$ as well as those cases
where eq.~(\ref{3.3}) has more than
one solution, assuming generic values for the free parameters.

In the second case where $Q$ contains a Jordanian 2-block,
equation (\ref{3.8}) in a first step reduces the set of independent
$E_{ijk}$ to the same 11 components as above.
Further analysis shows that $E_{111}$, $E_{133}$, $E_{211}$ and
$E_{311}$ vanish. The remaining components have to satisfy
\bea E_{333}&=&\b E_{333} \nonumber\\
     E_{222}&=&\a E_{222}+E_{122} \nonumber\\
     E_{322}&=&\a^2\b E_{322}+\a E_{123}+E_{132} \\
     E_{ijk}&=&\a_i\a_j\a_k E_{ijk} \q{, \ for} (ijk)=122,233,123,132\
     ,\nonumber\eea
with $\a_1=\a_2=\a$, $\a_3=\b$.
The solutions of this system constitute the cases 12 to 15 of the table.

In the third case, where $Q$ corresponds to a Jordanian 3-block,
5 of the 11 basic components vanish:
$E_{111}=E_{222}=E_{122}=E_{211}=E_{311}=0$. The remaining ones satisfy
\bea E_{333}&=&\a E_{333}+E_{233} \nonumber\\
    E_{233}&=&\a^3E_{233}+\a^2E_{133}+(\a+\a^3)E_{322}+\a E_{123} \nonumber\\
    E_{322}&=&\a^3 E_{322}+\a E_{123}+E_{132} \\
    E_{132}&=& - E_{123} \nonumber \\
    E_{ijk}&=&\a^3E_{ijk} \q{, \ for} (ijk)=133,123\ .\nonumber\eea
This gives the solutions 16 and 17.

In the table below an entry 0 in column $ijk$ means that $E_{ijk}$ has to
vanish in the respective case.
$(\a,\a,\b)^1$ denotes a matrix $Q$ of the second form in (\ref{3.9}) and
$(\a,\a,\a)^1_1$ of the third form.
$x,y\in\CC$, and $\g_n$ are arbitrary n-th roots of unity, i.e.~$(\g_n)^n = 1$.
The table is partially ordered, with the more special cases of $Q$
giving  less restrictions on the components of $E_{ijk}$
appearing first.

The proof that these quantum planes have the correct Poincar\'e series
to all degrees was done in \c{ArtSch}.
Those planes which play the role of quantum group modules will be discussed
in more detail in section 5.

\begin{center}\begin{tabular}{|c|c|c|c|c|c|c|c|c|c|c|c|}\hline
&$Q$     &\multicolumn{9}{|c}{$E_{ijk}$}&\raisebox{-.8ex}{123}\\
&$\a_i$  &111&222&333&211&311&122&322&133&233&132\\ \hline
1&$1,1,1$   &   &   &   &   &   &   &   &   &   &  \\ \hline
2&$1,1,\g_2$  &   &   &0  &   &0  &   &  0&   &   &0 \\ \hline
3&$1,\g_2,\g_2$ &   &0  &  0& 0 &0  &   &0  &   &0  &  \\ \hline
4&$1,\g_2,\g_4$ &   &0&  0&0  &0  &   &0  &0  &   &0   \\ \hline
5&$1,x,1/x$ &   & 0 &  0&0  &0  & 0 &0  &0  &0  &   \\ \hline
6&$\g_3,\g_3,\g_3$&0&0&0&   &   &   &   &   &   &   \\ \hline
7&$\g_6\ii^4,\g_6\ii^4,\g_6$&0&0&0&  &0  &   &0  &  &   &0  \\ \hline
8&$\g_9,\g_9\ii^7,\g_9\ii^4$& 0&0&0&   &0  & 0 &   &   &0  &0 \\ \hline
9&$x,x^{-2},x$&0  & 0 &0&   &0  & 0 &0  &0  &   &    \\ \hline
10&$x,x^{-2},\g_2 x$&0&0  &  0&  &0  & 0 &0  &0  &   &0   \\ \hline
11&$x,y,1/xy$ &0  & 0 &  0&0  &0  & 0 &0  &0  &0  &   \\ \hline
12&$(1,1,1)^1$&0&   &   &0&0&0  &   &0&   &$\S (5)$\\ \hline
13&$(\g_3,\g_3,\g_3)^1$&0&   &0  &0&0&$\S (1)$&   &0&   &$\S (5)$\\ \hline
14&$(\g_2,\g_2,1)^1$&0&0  &   &0&0&0  &  &0&0   &$\S (5)$\\ \hline
15&$(x,x,x^{-2})^1$&0&0  &0  &0&0&0  &  &0&0  &$\S (5)$\\ \hline
16&$(1,1,1)^1_1$&0&0& &0&0&0& &$\S (2)$&0  &$\S (5)$\\ \hline
17&$(\g_3,\g_3,\g_3)^1_1$&0&0& &0&0&0& &$\S (3)$&$\S (4)$&0\\ \hline
\multicolumn{12}{l}{
with ${\S (1)}{=}(1{-}\g_3)E_{222}\ ,\quad
      {\S (2)}{=}{-}E_{123}{-}2E_{322}\ ,\quad
      {\S (3)}{=}{-}(\g_3{+}\g_3\ii^2)E_{322}$}\\
\multicolumn{12}{l}{\phantom{with} ${\S (4)}{=}(1{-}\g_3)E_{333}\ ,\quad
{\S (5)}\ \a_1E_{123}{+}E_{132}{=}0 $.}
\end{tabular}\\[3mm]
{\bf Table 1: three-dimensional quantum planes}
\end{center}

\section{Quantum matrices in three dimensions}

Let $A^i_j$ be the generators of the matrix group acting on the quantum
plane $x^i$ from the left and preserving the plane in the sense that the
quantities $x'^i=A^i_jx^j$ satisfy again the quantum plane relations
(with $A^i_j$ and $x^k$ commuting).
Let $u_i$ be generating a coplane on which the quantum matrix acts
from the right by $u'_j=u_iA^i_j$. As for the plane we specify its
relations by
\be u_iu_j F_\a^{ij} =0 \qq{}(\ u\cd u\,F_\a=0\ ) \label{4.1}\ee
and demand a unique solution of the intersection equation
\be (\oo_i\ot F_\a)\, g^{i\a} = (F_\a\ot\oo_i)\, h^{\a i} = F\ ,\label{4.2}\ee
defining a tensor $(F^{ijk})$. Like for the plane we assume nondegeneracy
of $g$ and $h$ and get a cyclicity property for $F$
\be F^{ijk} = P^k_l F^{lij} \qq{,} P^k_l = h^{\a k}g_{\a l}\ .\label{4.3}\ee

As is well known, given a plane and a coplane in this way the relations for
the quantum group acting on them are fixed. Invariance of the relations
(\ref{3.1}) and (\ref{4.1}) gives
\be E^\a\, (A\cd A)\, x\cd x = 0 \qq, u\cd u\, (A\cd A)\, F_\a = 0 \
.\label{4.4}\ee
This should not lead to new relations for $x^i$ and $u_i$, i.e.
\bea E^\a\, A\cd A = D^\a_\b E^\b \qq, A\cd A\; F_\b  = F_\a \tilde D^\a_\b
\ .\label{4.6}\eea
The $D^\a_\b$ and $\tilde D^\a_\b$ will be called subdeterminants.

In the following we will consider the group algebra at increasing degrees:

{\bf i)  Degree 2.}
As a first implication from the Poincar\'e series of the group algebra
we show that the subdeterminants are independent and uniquely defined by
any one of eqs.~(\ref{4.6}). We proceed to show in detail how we
calculate the dimension of the algebra at a given order from the relations.
Eqs.~(\ref{4.6}) map a linear subspace of the homogeneous
component of order 2 of the algebra to new quantities $D^\a_\b$ and
$\tilde D^\a_\b$. This subspace, generated by the l.h.sides of (\ref{4.6})
is of dimension 45, because we have to subtract the 9 intersections of the
form $E^\a A\cd A\, F_\b$.
These intersections give rise to relations on $D^\a_\b$ and
$\tilde D^\a_\b$
\be D^\a_\g E^\g_{ij}F^{ij}_\b = E^\a_{ij}F^{ij}_\g \tilde D^\g_\b
\label{4.7}\ . \ee
To get the classical number of $\bino{10}2=45$ independent quantities
$A\cd A$, we need 9 of the quantities $D^\a_\b$ and $\tilde D^\a_\b$ to
be independent. This can only be accomplished with the matrix
$(E^\a_{ij}F^{ij}_\b)$ being invertible. We use the freedom of choosing
the upper Greek indices to take it to the unit matrix and have nine
independent subdeterminants $D^\a_\b = \tilde D^\a_\b$.
In this normalization we find from the contraction of (\ref{3.3}) with
(\ref{4.2}) the identity
\be e_{i\a}g^{i\a} = f_{\a i}h^{\a i}=E_{ijk}F^{ijk} =: \k \ .\label{4.8}\ee
We do not normalize $\k$ but choose a different normalization of $E$ and
$F$ further below.

{\bf ii)  Degree 3.}
The relations obtained from degree 2 are
\bea (E^\a\ot\oo)\,  A\cd A\cd A &=& D^\a_\b\cd A\, (E^\b\ot\oo)\label{4.10}\\
     (\oo\ot E^\a)\, A\cd A\cd A &=& A\cd D^\a_\b\, (\oo\ot E^\b)\label{4.11}\\
     A\cd A\cd A\, (F_\a\ot\oo)  &=& (F_\b\ot\oo)\, D_\a^\b\cd A\label{4.12}\\
     A\cd A\cd A\, (\oo\ot F_\a) &=& (\oo\ot F_\b)\,A\cd
     D_\a^\b\label{4.13}\ .\eea
The intersections of the l.h.sides of these equations are completely determined
by the intersection properties of plane and coplane separately, i.e.~by their
Poincar\'e series.
They give rise to relations between the quantities $D\cd A$ and $A\cd D$.
By comparison with (\ref{3.3}) we see that the intersection of the
first two equations must be given by the contraction of (\ref{4.10}) with
$f_{\a k}$ and of (\ref{4.11}) with $e_{i\a}$, leading to
\be f\, (D_\a\cd A_i) \, (E^\a\ot\oo^i) = e\, (A_i\cd D_\a) (\oo^i\ot E^\a)\ .
\label{4.14}\ee
Again by uniqueness of the solution of (\ref{3.3}) this gives
\be f\, (D\cd A)=\DD\,f \qq, e\, (A\cd D)  = \DD\, e \label{4.15}\ee
with $\DD$ a constant of proportionality. We will call $\DD$ the quantum
determinant of the group, which is hence given by
\be E\,(A\cd A\cd A) = \DD E\ .\label{4.16}\ee
The intersection of the other two equations similarly leads to a
proportionality constant $\tilde\DD$
\be (A\cd A\cd A)\, F=\tilde{\DD}F\ .\label{4.17}\ee
If $\k=EF\neq 0$ we conclude $\DD=\tilde\DD$.
Analogously to (\ref{4.15}) we find
\be (D\cd A)\, h = h\,\tilde{\DD} \qq, (A\cd D)\, g =g \tilde{\DD}\
.\label{4.23}\ee
Eqs.~(\ref{4.23}) together with (\ref{4.15}) mean that for
$\k\neq 0$ the following commutator relations hold
\be [D\cd A , h\ot f] = [A\cd D, g\ot e] = 0 \ .\label{4.23b}\ee
Two other relevant intersections are
$(E^\a\ot\oo)\, A\cd A\cd A\, (\oo\ot F_\b)$ and
$(\oo\ot E^\a)\, A\cd A\cd A\,(F_\b\ot\oo)$,
leading to
\bea (D\cd A)\, M=M\, (A\cd D) \qq, (A\cd D)\, N=N\, (D\cd A)\
,\label{4.25}\eea
with $M^{\a i}_{j\b}=E^\a_{jn}F^{ni}_\b$ and
$N^{i\a}_{\b j}=F_\b^{in}E_{nj}^\a$.
Eqs.~(\ref{4.25}) together imply
\be [MN,D\cd A] = [NM,A\cd D] = 0 \ .\label{4.26}\ee

Eqs.~(\ref{4.10})--(\ref{4.13}) map a linear subspace of $A\cd A\cd A$
of dimension $2\times 3^3\times 17 - {17}^2$ onto $D\cd A$ and $A\cd D$.
Hence the dimension in degree three is correct if the number of independent
$D\cd A$ and $A\cd D$ coincides with the classical number, which is 65.
Eqs.~(\ref{4.15}), (\ref{4.23}) and (\ref{4.25}) give all relations
of $A\cd D$ and $D\cd A$ following from intersections between
(\ref{4.10})--(\ref{4.13}).
Let us assume here for simplicity that they
contain complete ordering relations for $A\cd D$ and $D\cd A$.
The analysis should however be possible without this assumption.
Then the commutator relation (\ref{4.26}) has to leave more than 65 of the
$D\cd A$ independent. This restricts $MN$ to contain,
besides the unit matrix, at most a rank one operator.
By (\ref{4.23b}) this operator has to coincide with $h\ot f$ and hence
\be M N = \a\,\oo\ot\oo + \b\, h\ot f \label{4.27}\ .\ee
For $\k\neq 0$, where (\ref{4.23b}) is not valid, eqs.~(\ref{4.15})
and (\ref{4.23}) by themselves are too weak and we need $\b\neq 0$.
Then (\ref{4.26}) implies, also in this case, $\tilde\DD=\DD$, which thus
holds independently of $\k$.
The case $\b=0$ will be ruled out in the next section.

Analogous considerations for $A\cd D$ give
\be N M = \a\,\oo\ot\oo + \b\, g\ot e \label{4.28}\ ,\ee
with the same coefficients as in (\ref{4.27}). This follows from the
expression $MNM$ with the help of the following identities derived from
the cyclicity property of the planes
\be fM=e \q, Mg=h \q, eN=f \q, Nh=g \ .\label{4.29}\ee
With these identities, taking the trace in any of eqs.~(\ref{4.27}) or
(\ref{4.28}) gives \be \a+\b\k=1\ .\label{4.30}\ee

If $M$ and $N$ are nondegenerate, any of the eqs.~(\ref{4.25})
gives the complete commutation relations for $A\cd D$ and $D\cd A$
as we assumed above. One can show that to maintain this,
not both $M$ and $N$ must
be degenerate. Then with, say, $M$ invertible, eq.~(\ref{4.28}) with $\a=0$
gives $N=\k^{-1}g\ot f$. Another singular possibility for the r.h.s. of
(\ref{4.28}) is excluded because of the identities (\ref{4.29}).
Hence, excluding the singular cases $N{\sim}g\ot f$ and $M{\sim}h\ot e$ we
get that $M$ and $N$ must both be invertible as matrices in the composite
upper and lower indices and hence $\a\neq 0$.

Finally it is easy to verify that (\ref{4.27}) and (\ref{4.28}) imply the
existence of a Yang-Baxter operator.
The R-matrix with $A^{ij}_{kl}=F^{ij}_\a E^\a_{kl}$ as antisymmetrizer,
\be \hat R^{ij}_{kl} = \dd^i_k \dd^j_l - (1+q) F^{ij}_\a E^\a_{kl}
\ ,\label{4.31}  \ee
satisfies the Yang-Baxter equation iff $q$ satisfies
$q^2+q(2-\a^{-1})+1=0$.

{\bf iii)  Degree 4.}
Here we find the commutation relations of the
determinant with the group generators. The contractions
$(E\ot\oo)\, A\cd A\cd A\cd A \,(\oo\ot F)$ and
$(\oo\ot E)\, A\cd A\cd A\cd A \,(F\ot\oo)$ can be rewritten as
\be X\, A\cd\DD = \DD\cd A\, X \qq{and} Y\, \DD\cd A= A\cd\DD\, Y\
,\label{4.50}\ee
with the matrices $X$ and $Y$ given by
\be X^i_j = (E\ot\oo^i)(\oo_j\ot F) = e_{j\a}h^{\a i} \qq,
    Y^i_j = (\oo^i\ot E)(F\ot\oo_j) = g^{i\a}f_{\a j} \label{4.53}\ .\ee
$X$ and $Y$ obey with $P$ and $Q$ the identity
\be XQ = PY \label{4.54}\ .\ee
The quantum determinant is central if and only if $X$ is proportional to $\oo$.
Eqs.~(\ref{4.50}) both give the complete commutation relations
for the determinant, since $X$ and $Y$ are invertible matrices.
Since $\DD$ is assumed to be invertible in the group algebra and the $A^i_j$
independent this implies
\be XY = \l \oo\ .\label{4.51}\ee

A more systematic investigation of the degree 4 reveals more conditions on the
involved tensors. In particular it will turn out that the characteristic
matrices  $X$, $Y$, $P$ and $Q$ commute and that $X=Y^{-1}$ gives
an automorphism of the algebra. The relations in degree 4 are
\bea (E^\a\ot\oo\ot\oo)\, A\ot A\ot A\ot A &=&
     D^\a_\b\cd A\cd A\,(E^\b\ot\oo\ot\oo)\nn\\
     A\ot A\ot A\ot A\, (F_\b\ot\oo\ot\oo) &=&
     (F_\a\ot\oo\ot\oo)\, D^\a_\b\cd A\cd A
\label{4.32}\ ,\eea
together with the corresponding relations where $E^\a$ and $F_\b$ apply to
the other factors of $A$. Intersections between these relations we rewrite
with the help of the commutation relations for $D$ and $A$ as conditions on
$A\cd D\cd A$
\bea (e\ot\oo)\,A\cd D\cd A &=& \DD\cd A\,(e\ot\oo)   \label{4.33}\\
     (\oo\ot f)\, A\cd D\cd A &=& A\cd \DD\,(\oo\ot f) \label{4.34}\\
     E^\mu_{ip}F^{pq}_\a E^\nu_{qj}\,A^i_mD^\a_\b A^j_n &=&
     D^\mu_\sigma D^\nu_\tau\, E^\sigma_{mp}F^{pq}_\b E^\tau_{qn}
\label{4.35}\\
     A\cd D\cd A\, (g\ot\oo) &=& (g\ot\oo)\, \DD\cd A \label{4.36}\\
     A\cd D\cd A\, (\oo\ot h)&=& (\oo\ot h)\, A\cd \DD\label{4.37}\\
     A^i_mD^\a_\b A^j_n\, F^{mp}_\mu E^\b_{pq}F^{qn}_\nu &=&
F^{ip}_\sigma E^\a_{pq}F^{qj}_\tau \,D^\sigma_\mu D^\tau_\nu\ .\label{4.38}\eea

The relations for $D\cd D$, $\DD\cd A$ and $A\cd \DD$ following from
(\ref{4.33})--(\ref{4.38}) are:
\bea D^\a_\b D^\g_\dd \,E^\b_{mn}g^{n\dd} &=& E^\a_{ij}g^{j\g}\,A^i_m\DD
\label{4.45}\\
D^\a_\b D^\g_\dd\, h^{\b p}E^\dd_{pq} &=& h^{\a r}E^\g_{rl}\,\DD A^l_q
\label{4.46}\\
F^{ij}_\a e_{j\g}\,D^\a_\b D^\g_\dd &=& A^i_m\DD\, F^{mn}_\b e_{n\dd}
\label{4.47}\\
f_{\a k}F^{kl}_\g\, D^\a_\b D^\g_\dd &=& \DD A^l_q\,f_{\b p}F^{pq}_\dd
\label{4.48}\eea
together with the commutation relations for the determinant (\ref{4.50}).

Eq.~(\ref{4.45}) combined with (\ref{4.47}) gives similarly to
(\ref{4.51}), that the matrix
$U^i_j = e_{n\g}F^{in}_\a E^\a_{jm}g^{m\g}$ commutes with $A\cd\DD$ and
hence must be proportional to $\oo$.
The expression $E^\dd_{mi}F^{mj}_\epsilon$ gives a nondegenerate map from
indices $(i,j)$ to $(\dd,\epsilon)$, which is easily seen with the help of
the cyclicity conditions. Application to $U$ yields in the case
$\b\neq 0$ the new condition
$e_{j\g}h^{\dd j}f_{\epsilon n}g^{n\g}=\tilde\l \dd^\delta_\epsilon$, or
equivalently
\be YQ^{-1}XQ = \tilde\l\oo \ .\label{4.55}\ee
Together with eq.~(\ref{4.51}) we get $\l XQ=QX \tilde\l$.
For $\k=\hbox{tr}XQ\neq 0$ we thus find $\l=\tilde\l$ and
that $X$ and $Q$ commute. But then with (\ref{4.54}) and (\ref{4.51})
all characteristic matrices commute:
\be [P,Q]=[X,Q]=[Y,Q]=[X,P]=[Y,P]=0\ .\label{4.56}\ee
This property is the basis of our classification and
we will assume these commutation relations to hold even in the case
$\k=0$.

The dimension of the subspace of $A\cd A\cd A\cd A$ spanned by the l.h.sides of
(\ref{4.32}) is $2\times 3^4\times 66 - 66^2$, where 66 is the number of
relations for $x\cd x\cd x\cd x$.
Therefore the dimension in the quartics is correct if we have
the classical number of 270 independent quantities $A\cd D\cd A$.
The dimension of the subspace spanned by the l.h.sides of
(\ref{4.33})--(\ref{4.38}) turns out to be $2\times 3^3\times
12-12^2$ if we can prove that $e_{i\a}u_j$,
$v_if_{\a j}$ and $w_{\mu\nu}E^\mu_{ip}F^{pq}_\a E^\nu_{qj}$ span a space
of dimension 12, or equivalently that the equation
\be e_{i\a}u_j + v_if_{\a j} + w_{\mu\nu}E^\mu_{ip}F^{pq}_\a E^\nu_{qj}=0
\label{4.39}\ee
has exactly three solutions $(u,v,w)$.
This can be proven by multiplying (\ref{4.39}) by $N\ot\oo$ and defining
\be \tilde w_{\nu \b} = \a w_{\nu\b}-v_i F^{il}_\nu e_{l\b}\qq,
     \tilde u_j = u_j-\b w_{\a\b}h^{\a n}E^\b_{nj}\ .\label{4.41}\ee
Then (\ref{4.39}) becomes
\be f_{\nu k}\tilde u_j = \tilde w_{\nu\b}E^\b_{kj} \ .\label{4.42}\ee
By cyclicity $E^\a_{ij}\chi^j=0$ implies $\chi=0$.
Multiplying (\ref{4.42}) with $\chi^j$ and noting that $\tilde w$, if
nonzero, has to be nondegenerate (since $f$ is nondegenerate),
we find a contradiction in case $\tilde w$ and $\tilde u$ do not
vanish. Then with $\tilde w=\tilde u=0$, eqs.~(\ref{4.41}) give a
three-dimensional solution as we claimed.
We notice that $\b=0$ gives $u=0$, and  similarly by
multiplying (\ref{4.39}) with $1\ot M$, it gives $v=0$,
i.e.~(\ref{4.39}) has no solution.
Thus the case $\b=0$ has to be excluded.
We finally arrive at the correct number of independent $A\cd D\cd A$
if on the r.h.s. 45 of the $D\cd D$, $A\cd\DD$ and $\DD\cd A$ are
independent. This is again the classical number.

Eqs.~(\ref{4.45})--(\ref{4.48}) parallel, as relations for $D^\a_\b$,
the basic relations (\ref{4.6}) for $A^i_j$. We read
them as mapping some of the quantities $D\cd D$ onto the independent
quantities $\DD A$ and $A\DD$. Tensors like $E^\b_{mn}g^{n\dd}$ we
argued above to be independent for $m=1,2,3$.
As for $A\cd A$ we want 45 of the 81  $D\cd D$ to be independent.
The difference is that there now seem to be twice as many equations
compared with the case of $A\cd A$. To arrive at the correct number of
relations the l.h.sides of (\ref{4.45}) and (\ref{4.46}) therefore must
coincide, likewise for (\ref{4.47}) and (\ref{4.48}). In addition
$F^{ij}_\a e_{j\g}E^\a_{kn}g^{n\g}=\mu\dd^i_k$ must not vanish.
The proportionality between the l.h.sides is given through the matrices
$X$ and $Y$, with $\mu=\a\k+\b\lambda$:
\be E^\a_{mn}g^{n\b}Y^m_k=\mu h^{\a l}E^\b_{lk} \qq,
    F^{mn}_\a e_{n\b}X^k_m=\mu f_{\a l}F^{lk}_\b\ .\label{4.57}\ee
By lowering all indices of $E$ and raising all those of $F$ this turns into
\be E\,(Y\ot Y\ot Y) = \mu\lambda E \qq,
    (X\ot X\ot X)\, F =\mu\lambda F\ .\label{4.58}\ee
The antisymmetrizer which we introduced to write the R-matrix (\ref{4.31})
reads now
\be A^{ij}_{kl}:=F^{ij}_\a E^\a_{kl}=\l^{-1}E_{klm}X^m_nF^{nij}=
\l^{-1}F^{ijn}Y^m_nE_{mkl}\label{4.61}\ .\ee
$A$ satisfies $A^2=A$ and $\tr A=3$. From (\ref{4.51}) and (\ref{4.58})
we find $A (Y\otimes Y)=\mu^2\l^{-1}(Y\otimes Y)A$. Multiplying from left
and right with $A$, we obtain
\be [A,X\otimes X]=[A,Q\otimes Q]=0\ ,\label{4.62}\ee
together with $\l=\mu^2$.
The normalization of $X$ and $Y$, i.e.~of $g$ and $h$ relative to $e$ and
$f$ is still at our disposal.
{}From (\ref{4.57}) we see that $\mu$ is proportional to $X$ and $Y$.
We normalize
\be \mu=1 \qq{,\qquad and hence} \lambda=1 \q, \a=\b=\frac 1{1+\k} \
.\label{4.62b}\ee

Finally we rewrite the basic equation for degree 3 in terms of the
antisymmetrizer using the results of this section.
Contracting the free Greek indices in the first of eqs.~(\ref{4.28})
with $h_{a\a}$ and $f^{b\b}$ and using (\ref{4.51}) and
(\ref{4.58}), we get
\be (1+\k)\ A_{aj}^{ib}\ A^{cj}_{kd}\ (XQ)^{-1}\ii^d_b =
(XQ)^{-1}\ii^c_a\dd^i_k + \dd^i_a\dd^c_k \label{4.63}\ ,\ee
with $\k = E_{ijk}F^{ijk} = \tr XQ=\tr (XQ)^{-1}$.

Eq.~(\ref{4.55}), which together with (\ref{4.51}) implies the
commuting of all the characteristic matrix quantities, follows immediately
in degree 5. Starting from
$(E\ot\oo\ot\oo)\,A\cd A\cd A\cd A\cd A\,(\oo\ot\oo\ot F)$
we find
\be\DD D^\a_\b f_{\g p}g^{p\b}=f_{\dd k} g^{k\a}D^\dd_\g\DD\ .\label{4.97}\ee
and from $(\oo\ot\oo\ot E)\,A\cd A\cd A\cd A\cd A\,(F\ot\oo\ot\oo)$ we get
\be h^{\g k}e_{k\a}\DD D^\a_\b=D^\g_\dd\DD h^{\dd p}e_{p\b}\ .\label{4.98}\ee
With the independence of the $D^\a_\b$ eq.~(\ref{4.55}) follows
immediately. This derivation opens an easy way to generalize
(\ref{4.56}) to arbitrary dimensions, where $X$ is the contraction
of the $d$ dimensional $E$ and $F$ tensors over $d-1$ indices, and $Q$ and $P$
are their respective cyclicity matrices.

\section{Commutation relations of the determinant}

The conditions of the previous section allow transformations on the
solution space. This will be used to reduce the matrix $X$,
which gives the commutation relations of the quantum determinant
with the matrix elements, to one out of a few discrete cases.
For theses cases we will solve the conditions.
The transformations are induced by any automorphism of the algebra.
The possibility to find these without further knowledge of the solution
is due to eqs.~(\ref{4.58}) which determine $X$ to be
an automorphism of the algebra. This contains enough information
on the solutions to derive further automorphisms which give the desired
transformation on $X$.
In the end of section 4 we will derive the full automorphism group for each
solution, leading back to the general form of $X$.

Let $Z$ be an automorphism of the algebra, i.e.~an invertible mapping
$x^i\mapsto Z^i_jx^j$ and $u_i\mapsto u_j Z^{-1}\ii^j_i$,
leaving $E$ and $F$ invariant. Put
\be E'=E\,(Z^{-1}\otimes 1\otimes Z) \qq, F'=(Z\otimes 1\otimes Z^{-1})\,F\ .
\label{4.103}\ee
It is easy to see that the cyclicity conditions are satisfied with
\be Q' = Z^3 Q \qq, P' = Z^{-3} P \q{,\qquad and} [Q,Z] = [P,Z] =
0\label{4.104}\ ,\ee
the latter following from the uniqueness of $Q$.
{}From the definitions of $X$ and $Y$ we find
\be X' = Z^{-3} X \qq, Y'= Z^3 Y \q{,\qquad and} [Z,X] = [Z,Y] = 0
\label{4.107}\ .\ee
$(XQ)$ and $\k$ are invariant. With the new antisymmetrizer
\be A' = (1\otimes Z^{-1})\, A\, (1\otimes Z) \qq, [A,Z\otimes Z]=0\ ,
\label{4.109}\ee
it is now easy to verify that the transformation preserves all conditions of
the previous section.
In fact, transformation (\ref{4.103}) corresponds to a so
called twisting of a solution of the Yang-Baxter equation \c{Res}, with
\be \Rhat^F = \hat F\Rhat\hat F^{-1} \qq, F=Z\otimes 1 \q{, i.e.}
\hat F^{ij}_{kl}=Z^j_k \dd^i_l \ .\label{4.110}\ee
If $\Rhat=\oo-(1+q)A$ is a solution of the Yang-Baxter equation, so is
$\Rhat^F$, since $[\Rhat, Z\otimes Z] = 0$.

The matrix $X$ is an automorphism by virtue of eqs.~(\ref{4.58}) with the
normalization (\ref{4.62b}).
This condition allows one to classify $X$ in a way similar as it was done
for the quantum planes in section 1. We begin by assuming that $X$ is
diagonalizable.
Later we will prove that all other cases can be reduced to diagonal ones
by an appropriate twist transformation.
Then the condition involving the $E$ tensor becomes
\be E_{ijk} = \a_i\a_j\a_k E_{ijk}\ .\ee
As in section 1 this means that either a product of the eigenvalues of $X$
is equal to 1 or the corresponding component of $E$ has to vanish. This yields
table 2. Since eqs.~(\ref{4.58}) are invariant
under $X\mapsto\g_3X$, only one of the three cases is given.
The table is equally valid for the components of $F$.

The automorphism group of these algebras is determined by the same equations
which were just solved for $E$.
Restricting ourselves to automorphisms given by a diagonal matrix $Z$,
every nonvanishing component of $E$ or $F$  restricts a particular product
of eigenvalues of $Z$.
A subset of solutions is given by those $Z$ which coincide with $X$ up to
the choice of roots and parameters.
This permits in particular the twists given on the r.h.s. of the table below,
which transform the particular case to the one given in the last column
via $X'=Z^{-3}X$.
In the case $Z=X^{1/3}$ the roots of the eigenvalues have to be chosen
such that their triple products equals 1 if and only if this is the case for
the eigenvalues of $X$.

We see from table 2 that after making use of the twist transformation,
we have for $X$ the following seven cases of diagonal matrices:
\be X=(1,1,1),\ (\z_3,\z_3,\z_3),\ (1,\z_3,\z_3^2),\ (\z_9,\z_9^4,\z_9^7),\
      (1,1,\z_3),\ (\z_3,\z_3,1),\ (\z_3,\z_3,\z_3^2)\ ,\label{4.112}\ee
with $\z_n=e^{\pm 2\pi i/n}$.
We will show in the following that any nondiagonal $X$ can be brought into
diagonal form  by a twist transformation.
Hence (\ref{4.112}) is already the complete list of types of $X$-matrices.

Let us consider the case of $X$ containing a 2-block in Jordan normal form.
To solve simultaneously for the automorphisms, we consider eq.~(\ref{4.58})
with a slightly generalized $X$
\be   X=\bma{ccc}\a&\sigma&\\&\a&\\&&\b\ema \qq, \a,\b,\sigma\neq 0\ .
\label{4.113}\ee
Writing  (\ref{4.58}) in components for this $X$ one finds that (\ref{4.58})
is actually independent of the parameter $\sigma$.
Hence for any $X$ of the form (\ref{4.113}) with $\sigma=1$ there
exists the automorphism $Z$ of the same form and the same eigenvalues, but
with $\sigma=\frac 13$, transforming $X$ into
$X'=\hbox{diag}(\a^{-2},\a^{-2},\b^{-2})$.

In  case  $X$ consists of a Jordanian 3-block, we consider eq.~(\ref{4.58})
with
\be X=\bma{ccc}\a&\si&\tau\\&\a&\si\\&&\a\ema \qq, \a,\si\neq 0\ .
\label{4.115}\ee
The system of equations corresponding to (\ref{4.58}) now implies $\a^3=1$ and
depends on $\si$ and $\tau$ only through the combination
\be \g(X)= \dfrac{2\a\tau-\si^2}{\a\si}\ . \label{4.116}\ee
Defining for $n\in\QQ$ the  $n$-th power of the above $X$ to be
\be X^n=\a^n\bma{ccc}1&n\a^{-1}\si&\frac{n(n-1)}2\a^{-2}\si^2+n\a^{-1}\tau\\
0&1&n\a^{-1}\si\\0&0&1\ema\ ,\label{4.117}\ee
we find $\g(\lambda X^n)=\g(X)$, for $\lambda\in\CC$.
Hence if $X$ is an automorphism, so is $Z=\a^{-1/3} X^{1/3}$.
Twisting with this $Z$ transforms $X$ again to diagonal form.

\begin{center}\begin{tabular}{|c|c|c|c|c|c|c|c|c|c|c||c|c|}\hline
$X$&\multicolumn{9}{|c}{$E_{ijk}$,$F^{ijk}$}&\raisebox{-.8ex}{123}&
\multicolumn{2}{|c|}{twist}\\
$\a_i$              &111&222&333&211&311&122&322&133&233&132&$Z$&$X'$\\ \hline
$1,1,1$             &   &   &   &   &   &   &   &   &   &   &--&\\ \hline
$1,1,\g_3$          &   &   &   &   &0  &   &0  &0  &0  &0  &--&\\ \hline
$1,\g_3,\g_3^2$     &   &   &   &0  &0  &0  &0  &0  &0  &   &--&\\ \hline
$1,1,\g_2$          &   &   &0  &   &0  &   &0  &   &   &0  &$X$&1,1,1
\\ \hline
$1,\g_6^4,\g_6$     &   &   &0  &0  &0  &0  &0  &0  &   &0  &$X$&
$1,\g_3,\g_3^2$\\ \hline
$1,\g_2,\g_2$       &   &0  &0  &0  &0  &   &0  &   &0  &   &$X$&
1,1,1\\ \hline
$1,\g_4^2,\g_4$     &   &0  &0  &0  &0  &   &0  &0  &   &0  &$X^{-1}$&
1,1,1\\ \hline
$1,x^{-2},x$        &   &0  &0  &0  &0  &0  &0  &0  &   &0  &$X^{1/3}$&
1,1,1\\ \hline
$1,x^{-1},x$        &   &0  &0  &0  &0  &0  &0  &0  &0  &   &$X^{1/3}$&
1,1,1\\ \hline
$x^{-2},x,x$        &0  &0  &0  &0  &0  &   &0  &   &0  &   &$X^{1/3}$&
1,1,1\\ \hline
$x^{-2},x,\g_2x$    &0  &0  &0  &0  &0  &   &0  &   &0  &0  &$X$&
$x^{-2},x,x$\\ \hline
$\g_9^4,\g_9^7,\g_9$&0  &0  &0  &0  &   &   &0  &0  &   &0  &--&\\ \hline
$x^4,x^{-2},x$      &0  &0  &0  &0  &0  &   &0  &0  &   &0  &$X^{1/3}$&
1,1,1\\ \hline
$x,y,(xy)^{-1}$     &0  &0  &0  &0  &0  &0  &0  &0  &0  &   &$X^{1/3}$&
1,1,1\\ \hline
\end{tabular}\vspace{3mm}
{\bf Table 2: determinant commutation relations}
\end{center}

\newpage
\section{Solutions for GL(3)}

\begin{samepage}
So far we have classified quantum planes according to their characteristic
matrix $Q$. Furthermore, we have brought the determinant commutation relation
as determined by the matrix $X$ with the help of the twist transformation
to one out of seven possible diagonal forms.
Both classifications are compatible, since $Q$ and $X$ commute
and $X$ is diagonal. The simultaneous choice of $X$ and $Q$ is restricted
by the twist invariant condition $\tr(XQ) = \tr(XQ)^{-1}$, which fixes most
of the parameters in $Q$.
Given $X$ and $Q$, the characteristic matrix of the coplane follows
as $P=X^2Q$. From table 1 and 2 we can read off the independent components
of the tensors $E$ and $F$.

With these components we solve the remaining conditions (\ref{4.53}) and
(\ref{4.63}) for the matrix group which were derived in section 2 and
which we repeat here for convenience.
\bea E_{jmn}F^{mni} &=& X^i_j \label{A} \\
     (1+\k)\ A_{aj}^{ib}\ A^{cj}_{kd}\ (XQ)^{-1}\ii^d_b
     &=& (XQ)^{-1}\ii^c_a\dd^i_k + \dd^i_a\dd^c_k\label{B}\ ,\eea
where $\k=\tr XQ$ and $A^{ij}_{kl}=E_{klm}X^m_nF^{nij}$.
The second  eq.~(\ref{4.53}) follows from (\ref{A}) with $Y=X^{-1}$.

In some of the cases, where the matrices $X$ and $Q$ allow for nondiagonal
coordinate transformations, we introduce normal forms for the tensor
$E$, in order to solve eqs.~(\ref{A}) and (\ref{B}).
The solutions are further normalized  by rescaling the coordinates or using
more general coordinate transformations.

For every solution we will give explicitly the $E$ and $F$ tensors
from which the corresponding R-matrix is easily built.
According to section 2
\be \hat R^{ij}_{kl} = \dd^i_k \dd^j_l -(1+q)\,E_{klm}X^m_nF^{mij} \ee
satisfies the Yang-Baxter equation, with $q$ a solution of
\be q^2+q(1-\tr(XQ))+1=0\ ,\label{4.q}\ee
where the two solutions of this equation correspond to mutually inverse
R-matrices.

The results of our analysis are summarized in table 3 for
 $X$ and $Q$ which lead to quantum matrix groups.
There $(\a,\b,\g)$ denotes a diagonal matrix and $\z_n=e^{2\pi i/n}$.
For  $X=(\z_9,\z_9\ii^4,\z_9\ii^7)$ and $X=\g_3(1,1,\z_3)$  no solution
exists.
We recall that $X$ proportional to the unit matrix is equivalent to the
determinant being central. Out of the two solutions of
the quadratic eq.~(\ref{4.q}) we have listed only one.
$q=1$ represents the case of a double root and implies $\Rhat^2=1$.
Only the solutions B and F possess nontrivial parameters: solutions of type
B have 1 or 2 parameters (including $u$), solutions of type F have 1 parameter.
Cases where the twist invariant matrix $XQ$ is not diagonalizable will be
called Jordanian. F and G are Jordanian quantum groups, they all give
$\Rhat^2=1$.

As following from our construction,
the solutions listed in the table and given explicitly in the
following represent equivalence classes of solutions connected by twist
transformations. Since it is easy to determine the full automorphism group
of a given solution, the solutions derivable from the representatives
are also easy to obtain. We give the complete list of automorphisms at
the end of this section.
It may happen that the twist transformation introduces new parameters
into the solution.
The cases C, C' and E do not admit any further twist transformation.

\def\du{\rule[-1.7mm]{0cm}{5.5mm}}\def\dub{\rule[-7.5mm]{0cm}{17mm}}
\begin{center}\tabcolsep3mm\begin{tabular}{|c|c|c|c|}\hline
 &   $X$   &    $Q$     &$q$\\ \hline
A&         &\du$(1,1,1)$ & 1\\ \cline{1-1}\cline{3-4}
B&$(1,1,1)$&\du$(1,u,1/u)$&$u$\\ \cline{1-1}\cline{3-4}
C&         &\du$(\z_9,\z_9^4,\z_9^7)$&$\z_3$\\ \hline
C\makebox[0mm][l]{'}&\raisebox{-2.5mm}[0mm][0mm]{$(\z_3,\z_3,\z_3)$}&
\du$(\z_9,\z_9^4,\z_9^7)$&$\z_3$\\ \cline{1-1}\cline{3-4}
D&                  &\du$(1,\z_3,\z_3^2)$&$\z_3$\\ \hline
E&$(1,\z_3,\z_3^2)$&\du$(1,1,1)$&$\z_3$\\ \hline
F&\raisebox{-8mm}[0mm][0mm]{$(1,1,1)$}&
$\dub\bma{ccc}1&1&\\&1&\\&&1\ema$&1\\ \cline{1-1}\cline{3-4}
G&         &$\dub\bma{ccc}1&1&\\&1&1\\&&1\ema$&1\\ \hline
\end{tabular}\\[3mm]
{\bf Table 3: GL(3) quantum matrix groups}
\end{center}
\end{samepage}

\subsection{Cases with $X=\oo$}

Independent solutions are marked with a capital letter A...G and a number
which distinguishes subcases with the same matrices $X$ and $Q$.
Again we will denote a matrix in Jordanian normal form (\ref{3.9})
which contains a 2-block by $(\a,\a,\b)^\dpo$ and the Jordanian 3-block by
$(\a,\a,\a)^1_1$.

$X=\oo\q, Q=\oo$:\qquad
This being the case with the most number of independent variables,
we have at the same time the full GL(3) in the choice of
coordinates left.
By the cyclicity equation $E_{ijk}= E_{kij}$, the tensor $E$ decomposes into
\be E_{ijk} = \a^E\,\e_{ijk} + \phi^E_{ijk}\label{A.0}\ , \ee
where $\e$ is the usual completely antisymmetric tensor in 3 dimensions,
$\a^E$ a constant and $\phi^E$ a completely symmetric tensor.
We use the GL(3) freedom to take the completely symmetric tensor to one out of
the following 9 normal forms of ternary cubic forms \cite{Web}.
We give the corresponding monomial $\phi^E_{ijk}x^ix^jx^k$,
the verbal description refers to the corresponding curve
$\phi^Exxx=0$ in the projective plane:
\def\rtext#1{\makebox[0cm][l]{\makebox[12cm][r]{#1}}}
\bea
\phi^E&=&\rtext{(nondegenerate cubic)}x^3+y^3+z^3+6m\,xyz\q, 8m^3\neq-1
\hspace{4.5cm}\nn\\
\phi^E&=&\rtext{(cubic with self intersection)}x^3+y^3+6m\,xyz\q, m\neq 0
\nn\\
\phi^E&=&\rtext{(cubic curve with cusp)}x^3+y^2z\nn\\
\phi^E&=&\rtext{(quadratic curve and line, 2 intersection pts.)}
x^3+3xyz\nn\\
\phi^E&=&\rtext{(quadratic curve with tangential line)}xy^2+x^2z\nn\\
\phi^E&=&\rtext{(3 coinciding lines)}x^3\nn\\
\phi^E&=&\rtext{(3 lines, 2 coinciding)}x^2y\nn\\
\phi^E&=&\rtext{(3 lines, 3 intersection pts.)}xyz\nn\\
\phi^E&=&\rtext{(3 lines, 1 intersection pt.)}xy(x+y)\nn\eea
Besides this we have to consider the case $\phi^E=0$.
The cases consistent with conditions (\ref{A}) and (\ref{B})
are the following:

$\phi^E$ nondegenerate cubic ---
The equations exclude that all $F^{iii}$ vanish, i.e.~$F$ being just the
classical $\e$ tensor.
We have solutions related to a nondegenerate cubic curve, but they
turn out to be reducible to degenerate cases by twisting:
Besides $\det Z=1$ for $\a^E$ or $\a^F\neq 0$, $Z$ must be an automorphism
of the cubic curve, restricting it to permutations and multiplication with
$(1,\z_3,\z_3^2)$. Preservation of $X$ and $Q$ restricts the permutations
further to cyclic ones.
The corresponding twisting of the solutions leads either to
$m_E$ and $m_F$ corresponding to a degenerate cubic, or to the
factorized cubic $xyz$.
\\
$\phi^E$ quadratic curve and intersecting line ---
Part of the solution is contained in  (B2) below,
appearing for $Q=(1,u,1/u)$ at $u=1$. The solution not contained
there is
\beax{A1}
&E_{111}=1 \q, E_{123}=\dfrac{1-\z_3^2}3 \q, E_{132}=\dfrac{1-\z_3}3&\\
&F^{333}=\dfrac12\q, F^{123}=\dfrac{1-\z_3}2 \q, F^{132}=\dfrac{1-\z_3^2}2&
\eeax
$\phi^E$ three coinciding lines ---
Transformations in the 2-3 plane with unit determinant
leave the $E$-tensor invariant and affect only $F^{222}$, $F^{333}$,
$F^{233}$ and $F^{322}$. These four components represent a binary cubic
form which always factorizes and has 3 normal forms $x^3$, $x^2y$ and
$xy(x+y)$. We find solutions corresponding to the cubic form
$x^3$ or to the zero cubic form. The latter case however
is included in solution (B2), $Q=(1,u,1/u)$ at $u=1$.
Hence we stay here with the solution
\bex{A2}
E_{111}=E_{123}=-E_{132}=1\q, F^{222}=F^{123}=-F^{132}=\frac12
\eex
$\phi^E$ three lines, two coinciding ---
This gives two solutions
\beax{A3/A3}
&E_{211}=E_{123}=-E_{132}=1&\\
&F^{233}=F^{123}=-F^{132}=\frac12\q, F^{333}=0\q{or}\frac12&
\eeax
$\phi^E$ three lines, intersecting in three points ---
The solutions of this case are either contained in (B1) for
$Q=(1,u,1/u)$ at $u=1$, or coincide, after exchanging $E$ and $F$,
with (A1).
\\
$\phi^E=0$, $E$ classical $\e$-Tensor ---
Coplanes compatible with the classical $\e$-tensor we find after we
exchange the role of $E$ and $F$ and searching the solutions found for
$X=Q=\oo$ for such ones which allow the classical $\e$-tensor for $F$.
The only case we find is one which we have identified as a subcase
of (B2) at $u=1$, namely $E=\e+x^3$, $F=\e$.

$X=\oo\q, Q=(1,1,-1)$:\qquad
The components $E_{111}$, $E_{222}$, $E_{122}$ and $E_{211}$ describe a binary
cubic form which always factorizes. By a transformation in the 1-2 plane
we can take it into one of the three normal forms:
\bea
\phi^{E,12}&=&\rtext{(3 coinciding linear factors)}x^3\hspace{10.5cm}
\nn\\
\phi^{E,12}&=&\rtext{(3 linear factors, 2 coinciding)}x^2y\nn\\
\phi^{E,12}&=&\rtext{(3 indep. linear factors)}xy(x+y)\nn\eea
As a fourth possibility the four components $E_{111}$, $E_{222}$, $E_{122}$
and $E_{211}$ may vanish.
In all cases however the equations give a contradiction.

$X=\oo\q, Q=(1,u,1/u)$:\qquad
$u=-1$ corresponds to the case $Q=(1,-1,-1)$ below and we exclude it here.
Then we stay with two solutions
\bex{B1}
E_{123}=1\q, F^{123}=\frac1{u+1} \q, E_{132}F^{132}=\frac{u}{u+1}
\eex
\bex{B2}
E_{123}=E_{111}=1\q, F^{123}=\frac1{u+1}\q,
E_{132}=-u^{1/3}\q, F^{132}=-\frac{u^{2/3}}{u+1}
\eex
Solution (B1) has in addition to $u$ another nontrivial parameter
$\nu=-E_{132}$.

$X=\oo\q, Q=(\z_9,\z_9\protect\ii^7,\z_9\protect\ii^4)$:\qquad
We find the solution
\beax{C1}
&&E_{133}=E_{322}=E_{211}=1 \q, F^{133}=\si \q, F^{322}=\si-\z_9\ii^4 \q,
F^{211}=\si+\z_9 \\
&& \hbox{with } \si=\frac13(-\z_9\ii^6+\z_9\ii^4+\z_9\ii^3-\z_9)
\eeax

$X=\oo\q, Q=(1,1,1)^\dpo$:\qquad
We make use of coordinate transformations commuting with $X$ and $Q$
to reduce the number of parameters in the solutions. We use
rescaling with $(a,a,a^{-2})$ and transformation with
$\bsma{ccc}1&&a\\&1&\\&&1\esma$. Then we stay with the solutions
\beax{F1-F4}
&E_{123}=-E_{132}=1\q, E_{222}=0\q{or}1\q, E_{322}=-\frac5{12}\mp\frac1{12}\\
&F^{123}=-F^{132}=\frac12\q, F^{311}=\frac16 \pm \frac1{12}\q,
F^{333}=\nu\neq0& \eeax
\beax{F5-F8}
&E_{123}=-E_{132}=1\q, E_{222}=0\q{or}1\q, E_{322}=\nu \hbox{ (arb.)}&\\
&F^{123}=-F^{132}=\frac12\q, F^{311}=\frac12 E_{322}+\frac38\pm\frac18&
\eeax
\beax{F9-F12}
&E_{123}=-E_{132}=1\q, E_{222}=0\q{or}1\q, E_{322}=-\frac34\pm\frac14& \\
&F^{123}=-F^{132}=\frac12\q, F^{111}=\nu \hbox{ (arb.)}&
\eeax

$X=\oo\q, Q=(1,1,1)^1_1$:\qquad
By rescaling and transforming with
$\bsma{ccc} 1&a&b\\&1&a\\&&1\esma$, which commutes with $Q$, we remove
all free parameters and have solutions
\beax{G1}
&E_{123}=-E_{132}=1\q, E_{333}=-\frac1{27}\q, E_{133}=-\frac13\q,
E_{322}=-\frac13&\\
&F^{123}=-F^{132}=\frac12\q, F^{311}=\frac13\q, F^{122}=\frac1{12}&
\eeax
\beax{G2}
&E_{123}=-E_{132}=1\q, E_{133}=-1&\\
&F^{123}=-F^{132}=\frac12\q, F^{311}=\frac12&
\eeax

$X=\oo\q, Q=(1,-1,-1)\q, Q=(\z_4,\z_4\protect\ii^2,\z_4\protect\ii^3)\q{or}
Q=(-1,-1,1)^\dpo$:\qquad
These remaining cases of $Q$ in table 1, allowed by the condition
$\tr XQ=\tr(XQ)^{-1}$, all lead to contradicting equations.

\subsection{Cases with $X=\z_3\oo$}

$X=\z_3\oo\q, Q=(1,\z_3,\z_3\protect\ii^2)$:\qquad
Part of the solutions is contained in the case
$Q=(u,\z_3\ii^2,\z_3/u)$ at $u=1$, see below.
Here remains
\beax{D1/D2}
&E_{123}=1\qq, E_{132}=-\g_3\qq, E_{111}=1&\\
&F^{123}=-\z_3\qq, F^{132}=\g_3\ii^2\qq, F^{222}=0\q{or} 1&
\eeax

$X=\z_3\oo\q, Q=\z_3\protect\ii^2\oo$:\qquad
This case leaves the full GL(3) for transformations of the coordinates.
To arrive at a normal form for $E$ we make use of the decomposition
(\ref{4.101a}) below. Imposing the cyclicity equation reduces $E$ and $F$
to the form
\be E_{ijk}=S^m_i\e_{mjk} -\z_3\e_{ijn}S^n_k\qq,
    F^{ijk}=\t S_m^i\e^{mjk} -\z_3\e^{ijn}\t S_n^k\ ,\ee
with $\tr S=\tr \t S=0$.
Condition (\ref{A}) then reads $[S,\t S]=0$ with $\tr(S\t S)=-\z_3\ii^2$.
$S$ and $\t S$ transform under a change of coordinates as matrices and we
can take $S$ into one of 3 cases of Jordan normal form:

$S=(\mu,\nu,-\mu-\nu)$ --- We may assume $\t S$ to be diagonalizable $\t
S=(\t\mu,\t\nu,-\t\mu-\t\nu)$, since otherwise we exchange $E$ and $F$.
The resulting solution is contained
in the solution for $Q=(u,\z_3\ii^2,\z_3/u)$ with $u=\z_3\ii^2$,
see below.
\\
$S=(\mu,\mu,-2\mu)^\dpo$ --- We take $\t S$ arbitrary.
The only solution found admits automorphisms $Z=(\a,\a,\a^{-2})$.
Twisting with
$\a=\z_9$ ($\z_9\ii^3=\z_3$) transforms $X$ and $Q$ into
$X=\oo$, $Q=\oo$. ---
$S=(0,0,0)^1_1$ gives no solutions.

$X=\z_3\oo\q, Q=(\z_9,\z_9\protect\ii^7,\z_9\protect\ii^4)$:\qquad
This case is completely parallel to the case $X=\oo$,
$Q=(\z_9,\z_9\ii^4,\z_9\ii^7)$. Taking $\z_3=\z_9\ii^3$, the solution is
\beax{C2}
&&E_{133}=E_{211}=E_{322}=1 \q, F^{133}=\si \q, F^{322}=\si-\z_9\ii^7 \q,
F^{211}=\si+\z_4 \\
&& \hbox{with } \si=\frac13(\z_9\ii^7-\z_9\ii^6-\z_9\ii^4+\z_9^3)
\eeax

$X=\z_3\oo\q, Q=(u,\z_3\protect\ii^2,\z_3/u)$:\qquad
The corresponding solution allows for automorphisms $Z=(\a,\b,1/\a\b)$.
Taking $\a=\b=\z_9$ ($\z_9\ii^3=\z_3$), we can transform $X$ and $Q$ to
$X=\oo$, $Q=(\z_3u,1,\z_3\ii^2/u)$.

$X=\z_3\oo\q,Q=(\z_3\protect\ii^2,\z_3\protect\ii^2,\z_3\protect\ii^2)^\dpo$:
\qquad
A coordinate transformation with
\raisebox{2ex}[2ex][0pt]{$\bsma{ccc}1&&a\\&1&\\&b&1\esma$},
commuting with $Q$, allows us to get rid of the components $E_{222}$,
$E_{233}$.
The solutions then allow automorphisms with $Z=(\a,\a,\a^{-2})$.
Taking $\a=\z_9$ ($\z_9\ii^3=\z_3$), we can twist $X$ and $Q$ into
$X=\oo$ and $Q$ a matrix with Jordan 2-block and eigenvalues 1.
Hence this solution is reducible to one of the cases $X=\oo$.

$X=\z_3\oo\q,
Q=\z_3\protect\ii^2(1,1,-1)$ or
$\z_3\protect\ii^2(-1,1,-1)$ or
$(\z_{12},\z_{12}\protect\ii^{10},\z_{12}\protect\ii^7)$ or
$(-\z_3\protect\ii^2,-\z_3\protect\ii^2,\z_3\protect\ii^2)^\dpo$ or
$(\z_3\protect\ii^2,\z_3\protect\ii^2,\z_3\protect\ii^2)^1_1$:\qquad
The first case of $Q$ leads to contradicting equations while in the other
cases there are not enough components to allow for a solution
of eq.~(\ref{A}).

\subsection{Cases with noncentral determinant}

$X=(1,\z_3,\z_3\protect\ii^2)$ reduces the $E$ tensor already to components
$E_{iii}$, $E_{123}$ and $E_{132}$.
$Q$ must be diagonal in order to commute with $X$.

$X=(1,\z_3,\z_3\protect\ii^2) \q, Q=\oo$:\quad
Part of the solutions are contained in the case
$Q=(1,u,1/u)$ at $u=1$ below.
Remaining here are the solutions
\beax{E1}
&E_{123}=1\q, E_{132}=-\z_3\ii^{1/3}\q, E_{222}=E_{333}=1&\\
&F^{123}=-\z_3\q, F^{132}=\z_3\ii^{5/3}\q, F^{111}=\tfrac13(1-\z_3\ii^2)&
\eeax
\beax{E2}
&E_{123}=1\q, E_{132}=-\z_3\ii^{1/3}\q, E_{222}=1&\\
&F^{123}=-\z_3\q, F^{132}=\z_3\ii^{5/3}&
\eeax

$X=(1,\z_3,\z_3\protect\ii^2)\q, Q=(1,-1,-1)$:\qquad
With $P=(1,-\z_3\ii^2,-\z_3)$ we get the same independent components for
the tensors $E$ and $F$ as in
the case $Q=(1,u,1/u)$. Therefore all solutions will be contained there
as the special case $u=-1$.

$X=(1,\z_3,\z_3\protect\ii^2)\q, Q=(1,u,1/u)$:\qquad
The two solutions allow automorphisms with $Z=(\a,\b,1/\a\b)$,
in one case restricted to $\a^3=1$. Taking $\a=1$, $\b=\z_9$
($\z_9\ii^3=\z_3$),
transforms $X$ and $Q$ into $X=\oo$, $Q=(1,\z_3u,\z_3\ii^2/u)$.

$X=(1,\z_3,\z_3\protect\ii^2)\q, Q=(\z_3,\z_3,\z_3)$:\qquad
Some solutions are included in the case $Q=(u,\z_3\ii^2/u,\z_3)$ with
$u=\z_3$. The remaining solution allows automorphisms
$Z=(\a,\g_3,\g_3\ii^2/\a)$.
Taking $\a=\z_9\ii^2$ ($\z_9\ii^3=\z_3$) transforms $X$ and $Q$ into
$X=\z_3\oo$, $Q=(1,\z_3,\z_3\ii^2)$.

$X=(1,\z_3,\z_3\protect\ii^2)\q, Q=(-\z_3,-\z_3,\z_3)$:\qquad
With $P=(-\z_3,-1,\z_3\ii^2)$ we get the same independent components
for the tensors $E$ and $F$ as in
the case $Q=(u,\z_3\ii^2/u,\z_3)$. Therefore all solutions will be
contained there as the special case $u=-\z_3$.

$X=(1,\z_3,\z_3\protect\ii^2)\q, Q=(u,\z_3\protect\ii^2/u,\z_3)$:\qquad
Depending on the value of $u$ we have two different solutions.
For $u=1$ we find back a solution of the case
$Q=(1,u,1/u)$ with $u=\z_3\ii^2$. For the other solution we need $u\neq
1,-1$. This solution however allows automorphisms $Z=(\a,\b,1/\a\b)$.
Taking $\a=1$, $\b=\z_9$, $\z_9\ii^3=\z_3$,
we twist $X$ and $Q=(u,\z_3\ii^2/u,\z_3)$ into $X=\oo$, $Q=(u,1/u,1)$.

$X=(\z_9,\z_9\protect\ii^4,\z_9\protect\ii^7)$:\qquad
In all cases the equations give a contradiction.

$X=\g_3(1,1,\z_3)$:\qquad
This $X$ restricts the independent components of $E$ to $E_{iii}$,
$E_{122}$ and $E_{211}$. However we do not find any $Q$
which leaves enough components of $E$ free and satisfies at the same time
$\tr XQ=\tr(XQ)^{-1}$. Especially $Q=\g_3\ii^2(1,1,\z_3\ii^2)$
leads to an $E$ with a non unique solution for the cyclicity equation.

\subsection{Automorphisms and twists}

For arbitrary $E_{ijk}$ and $F^{ijk}$ we have unique decompositions
\bea E_{ijk}&=&\a\,\e_{ijk} + S^m_i\e_{mjk} + \e_{ijn}T^n_k+\phi_{ijk}
\label{4.101a}\\
     F^{ijk}&=&\t\a\,\e^{ijk}+\t S_m^i\e^{mjk}+\e^{ijn}\t T_n^k+\t\phi^{ijk}
\label{4.101b}\ ,\eea
with $\e$ the usual $\e$-tensor, $S$, $\t S$, $T$ and $\t T$ traceless
matrices and $\phi$, $\t\phi$ completely symmetric tensors.
$Z$ being an automorphism is then equivalent to
\bea &z\,SZ=ZS\q, z\,ZT=TZ\q, z\,Z\t S=\t SZ\q, z\,\t T Z=Z\t T&\nn\\
     &\phi(Z\otimes Z\otimes Z)=\phi\qq, (Z\otimes Z\otimes Z)\t\phi=\t\phi&
\label{4.101c}\eea
where $z=\det Z$. If $\a$ or $\t\a\neq 0$ we have $\det Z=1$.
Eqs.~(\ref{4.101c}) imply $\det Z=1$ also if any of the
determinants of $S$, $T$, $\t S$ or $\t T$ does not vanish.

For all solutions found so far the automorphism group is given below.
The automorphisms of the cubic forms $\phi$, $\t\phi$ can be found in the
literature \c{Web}.
$Z$ is determined only up to a factor $\g_3$ which we have omitted.
Since any automorphism induces a twist (\ref{4.103}),
we determine here at the same time the full space of solutions given by twists
of the standardized solutions above.

\newpage
\begin{samepage}
\begin{tabular}{lll}
(A1)&$Z=(1,\g_3,\g_3\ii^2)$&\\[1mm]
(A2)&$Z=\bma{ccc} 1&&\\z_{21}&\g_3&z_{23}\\z_{31}&&\g_3\ii^2\ema$&\\[6mm]
(A3)&$Z=\bma{ccc} x&&\\&x^{-2}&\\z_{31}&&x\ema$&\\[6mm]
(A4)&$Z=\bma{ccc} x&&\\&x^{-2}&\\z_{31}&\frac{1-x^3}{3x^2}&x\ema$&\\[6mm]
(B1)&$Z=\Pi\circ (x,y,1/xy)$&
$\begin{array}{ll}u\neq 1:&\Pi=\oo\\[1mm]
                  u=1,\ \nu=-1:&\Pi=\hbox{arb. permutation}\\[1mm]
                  u=1,\ \nu\neq 1,-1:&\Pi=\hbox{cyclic perm.}\\[1mm]
                  u=\nu=1:&Z\in SL(3) \end{array}$ \\[10mm]
(B2)&$Z=(1,x,x^{-1})$&$\begin{array}{ll}u^{1/3}=1:&Z\in SL(3)\end{array}$
\\[1mm]
(C1), (C2)&--&\\[1mm]
(D1)&$Z=(1,x,x^{-1})$&\\[1mm]
(D2), (E1), (E2)&$Z=(1,\g_3,\g_3\ii^2)$&\\[1mm]
(F1-4)&$Z=\bma{ccc} \pm1&z_{12}&\\&\pm1&\\&z_{32}&1\ema$&
$\begin{array}{ll}\hbox{(F1):}&z_{32}=0\\[1mm]
                  \hbox{(F2):}&z_{32}=-2\pm 2\\[1mm]
                  \hbox{(F4):}&\hbox{only $+$ signs} \end{array}$ \\[8.5mm]
(F5-8)&$Z=\bma{ccc} x&z_{12}&z_{13}\\&x&\\&z_{32}&x^{-2} \ema$&
$\begin{array}{lll}\hbox{(F5):}&\nu\neq-\frac23:&z_{13}=0\\[1mm]
                              &\nu\neq-\frac13:&z_{32}=0\\[1mm]
                   \hbox{(F6):}&\nu\neq-\frac23:&z_{13}=0\\[1mm]
                     &\nu\neq-\frac13:&z_{32}=\frac{1-x^3}{x^2(3\nu+1)}\\[1mm]
                              &\nu=-\frac13:&x=1\\[1mm]
                   \hbox{(F7):}&\nu\neq-\frac16:&z_{13}=0\\[1mm]
                              &\nu\neq-\frac13:&z_{32}=0\\[1mm]
                   \hbox{(F8):}&\nu\neq-\frac16:&z_{13}=0\\[1mm]
                     &\nu\neq-\frac13:&z_{32}=\frac{1-x^3}{x^2(3\nu+1)}\\[1mm]
                              &\nu=-\frac13:&x=1 \end{array}$\\[29mm]
(F9-12)&$Z=
\bma{ccc} x&z_{12}&\frac{2\nu(x^3-1)}{x^2}\\&x&\\&z_{32}&x^{-2}\ema$&
$\begin{array}{ll}\hbox{(F9), (F11):}&z_{32}=0\\[1mm]
                  \hbox{(F10):}&z_{32}=2(x-x^{-2})\\[1mm]
                  \hbox{(F12):}&z_{32}=(x-x^{-2})/2 \end{array}$\\[6mm]
(G1), (G2)&$Z=\bma{ccc} 1&x&\frac12x(x-1)\\&1&x\\&&1\ema$&
\end{tabular}
\begin{center}
{\bf Table 4: automorphisms of the GL(3) solutions}
\end{center}
\end{samepage}

\newpage
\section{Normal ordering and Poincare series}

All quantum groups found here, with exception of the cases C,C', allow for an
alphabetic ordering of monomials. This provides a simple possibility to prove
the correctness of the full Poincare series with the help of the
diamond lemma \c{Ber}. To this end we have to show that ordering triple
products of generators always leads to a unique result.

Since it turns out that we can construct an ordering for the group generators
from that of the planes, we first treat the corresponding planes.
For a proof of the correctness of the Poincar\'e series of the planes in the
cases C,C' by different methods we refer to \c{ArtSch}.
The twisting transformation preserves the Poincare series of a solution,
as we will show in the end of this section.

Given a quantum plane and an ordering for the generators such as
\bea x^1<x^2<x^3\q, x_3<x_2<x_1 && \hbox{(solutions A, B, D)}\nn\\
     x^3<x^2<x^1\q, x_1<x_2<x_3 && \hbox{(solutions E, G)}\nn\\
     x^2<x^1<x^3\q, x_3<x_1<x_2 && \hbox{(solution F)}\eea
we can rewrite the relations $E^\a\,x\cd x$ as ordering rules by solving
for the lexically highest monomial in one of the relations,
\be x^u\cd x^v = \t E^\a\,x\cd x \qq,
\t E^\a = \oo^u\ot\oo^v - (E^\a_{uv})^{-1} E^\a \qq,\a=1\ldots 3\ ,\ee
eliminating it from the remaining relations.
If there exists some ordering of the generators where we get in this way
a system which allows for the replacement of the three anti-ordered
products, e.g. $x^3x^2$, $x^3x^1$ and $x^2x^1$ in cases A, B and D, we
say that the relations allow for an alphabetic ordering.
The orderings of generators given above are just chosen that this is the case.
In case C one can show that this is not possible, i.e.~there exist no
nondegenerate matrices $m$ and $n$ such that the three relations
$m^\a_\g E^\g_{kl}$ with an ordering $n^1_ix^i<n^2_jx^j<n^3_kx^k$
lead to appropriate substitution rules.

By the Diamond Lemma, a quantum plane which allows for an alphabetical
ordering has the classical Poincar\'e series if the ambiguously orderable
triple products, e.g. in cases A, B and~D $x^3x^2x^1$
which can be ordered either by starting with
$x^3x^2$ or with $x^2x^1$, corresponds to a unique ordered expression.
This we verified for all respective planes and coplanes to our solutions
for the above orderings. We observe that the coplane of a
quantum group always accepts the ordering opposite to that of the plane.

For the group generators we have the relations
\be E^\a\,(A\cd A)\, S_{ij} = 0 \qq, S^{ij}\, (A\cd A)\, F_\a=0 \qq{for}
\a,i,j=1,2,3\ ,\label{5.1}\ee
where $S^{ij}_{kl}=\dd^i_k\dd^j_l - A^{ij}_{kl}=\dd^i_k\dd^j_l -
E_{kla}X^a_bF^{bij}$ is the symmetrizer. We use the following orderings:
\begin{samepage}
\bea
A^1_3<A^2_3<A^3_3<A^1_2<A^2_2<A^3_2<A^1_1<A^2_1<A^3_1 &&
\hbox{(solutions A, B, D)}\nn\\
A^3_1<A^2_1<A^1_1<A^3_2<A^2_2<A^1_2<A^3_3<A^2_3<A^1_3 &&
\hbox{(solutions E, G)}\nn\\
A^2_1<A^1_1<A^3_1<A^2_3<A^1_3<A^3_3<A^2_2<A^1_2<A^3_2 &&
\hbox{(solution F)}\eea
\end{samepage}
The ordering of the upper indices coincides here with that of the plane
generators, for the lower indices with the coplane ones,
having given higher lexical significance to the lower indices.
Like for the planes we can rewrite relations (\ref{5.1}) into a
substitution system.
By inspection it is proven that the above orderings just lead to substitution
rules for the 36 anti-ordered squares. Finally we have to prove for
$\bino93=84$ anti-ordered triple products of generators that ordering with this
substitution system leads to a unique result. Then we conclude with the Diamond
Lemma that the Poincar\'e series of the algebra of group generators
coincides with that of 9 commuting variables.
This can easily be done on a computer.

Let us consider the effect of twisting on the Poincar\'e series of the planes.
Twisting with $Z$ changes $E^\a$ into
\be E^\a \mapsto z^\a_\b E^\b(1\otimes Z)\qq, z^\a_\b=e^{\a a} Z^{-1}\ii^b_a
e_{b\b} \ .\ee
The condition that $Z$ is an automorphism reads with $E^\a$:
\be E^\a(Z\otimes Z) = z^\a_\b E^\b\ ,\label{5.5}\ee
similar relations hold for $F_\b$.
The Poincar\'e series of the plane is determined by the linear dependences
between the tensors $(1\otimes\ldots \otimes E^\a\otimes \ldots\otimes 1)$,
i.e.~in the $n$-th order of the algebra, by the number of independent
solutions of the equation
\be E^\a_{i_1i_2}e^{(1)}_{\a i_3\ldots i_n} +
E^\a_{i_2i_3}e^{(2)}_{i_1\a i_4\ldots i_n} +
E^\a_{i_{n-1}i_n}e^{(n-1)}_{i_1\ldots i_{n-2}\a}\ .\label{5.6}\ee
For the twisted tensors $E^\a(1\otimes Z)$ this equation becomes
\be E^\a_{i_1j}Z^j_{i2}e^{(1)}_{\a i_3\ldots i_n} +
    E^\a_{i_2j}Z^j_{i_3}e^{(2)}_{i_1\a i_4\ldots i_n} +
    E^\a_{i_{n-1}j}Z^j_{i_n}e^{(n-1)}_{i_1\ldots i_{n-2}\a} \ .\label{5.7}\ee
Multiplying the last equation with
$(Z^{n-1}\otimes\ldots\otimes Z^1 \otimes 1)$
and using (\ref{5.5}) we see that the solutions of (\ref{5.6}) and
(\ref{5.7}) are in  1--1 correspondence and the Poincar\'e series of twisted
and untwisted plane coincide.

The relations of the group have the form
$(1\otimes..\otimes E^\a\otimes..\otimes 1) \otimes
(1\otimes..\otimes F_\b\otimes..\otimes 1)$.
The Poincar\'e series of the group is determined by linear dependences
between relations with different positions of the tensors $E$ and $F$ and
different values of $\a$ and $\b$.
These linear dependences however are like in the case of the planes not
changed by twisting, since we can reduce
the effect of twisting by the same manipulations as above
to the multiplication of the different factors
in the tensor product by increasing powers of $Z$, which can be removed by
redefinition. Therefore the Poincar\'e series of the group as well is not
changed by twisting.

\newpage
\appendix

\section{Appendix: R-Matrices}

The following R-matrices are derived from the solutions in the section 4 by
$\hat R^{ij}_{kl} = \dd^i_k \dd^j_l -(1+q)\,E_{klm}X^m_nF^{mij}$,
with $q$ a solution of $q^2+q(1-\tr(XQ))+1=0$.
In the cases $\tr(XQ)=\ 3,\ 0,\ 1+u+u^{-1}$ we took
$q=\ 1,\ \z_3\ii^2,\ u$.
Not all R-matrices are given here.

\vfill

(A1)\\*[-6ex]\ii
\benon\arraycolsep3mm\bma{ccccccccc}
   1   &       &       &       &       &       &       &       &       \\
       &\z_3^2 &       &       &       &       &       &       &       \\
       &       &\z_3   &       &       &       &       &       &       \\
       &       &       &\z_3   &       &       &       &       &       \\
       &       &       &       &   1   &       &       &       &       \\
\o{\ \z_3^2{-}1}&      &       &       &       &\z_3^2 &       &       & \\
       &       &       &       &       &       &\z_3^2 &       &         \\
\o{\ \z_3{-}1}&       &       &       &       &       &       &\z_3   &  \\
       &\o{\tfrac13(\z_3^2{-}1)}&       &\o{\tfrac13(\z_3{-}1)}& & & & & 1
\ema\eenon
\\[2ex]
\partcolumn{7.5cm}{
(A2)\\*[-2ex]\ii
\benon\arraycolsep2mm\bma{ccccccccc}
   1   &       &       &       &       &       &       &       &       \\
       &   1   &       &       &       &       &       &       &       \\
       &       &   1   &       &       &       &       &       &       \\
       &       &       &   1   &       &       &       &       &       \\
       &       &   1   &       &   1   &       & \o{{-}1}&       &       \\
1      &       &       &       &       &   1   &       &       &       \\
       &       &       &       &       &       &   1   &       &       \\
\o{{-}1} &       &       &       &       &       &       &   1   &       \\
       &       &       &       &       &       &       &       &   1
\ema\eenon}\restcolumn{
(A3/A4)\\*[-2ex]\ii
\benon\arraycolsep2mm\bma{ccccccccc}
   1   &       &       &       &       &       &       &       &       \\
       &   1   &       &       &       &       &       &       &       \\
 \o{{-}1}&       &   1   &       &       &       &       &       &       \\
       &       &       &   1   &       &       &       &       &       \\
       &       &       &       &   1   &       &       &       &       \\
       &       &       &   2   &       &   1   &       &       &       \\
   1   &       &       &       &       &       &   1   &       &       \\
       & \o{{-}2}&       &       &       &       &       &   1   &       \\
 \o{{-}1}& \o{{-}\si}&   1   & \si &       &       & \o{{-}1}&       &   1
 \ema\eenon \rightline{$\si=0$ (A3) or 1 (A4)}
}
\\[3ex]
(B1)\\*[-6ex]\ii
\benon\bma{ccccccccc}
   1   &       &       &       &       &       &       &       &       \\
       &u/\nu  &       &       &       &       &       &       &       \\
       &       &  \nu  &       &       &       &\o{1{-}u}  &       &       \\
       &\o{1{-}u}  &       &  \nu  &       &       &       &       &       \\
       &       &       &       &   1   &       &       &       &       \\
       &       &       &       &       &u/\nu  &       &\o{1{-}u}  &       \\
       &       &       &       &       &       &u/\nu  &       &       \\
       &       &       &       &       &       &       & \nu   &       \\
&&\phantom{u/\nu}&\phantom{u/\nu}&\phantom{u/\nu}&&&\phantom{u/\nu}&1
\ema\eenon

\newpage

(B2)\\*[-6ex]\ii
\benon\bma{ccccccccc}
   1   &       &       &       &       &       &       &       &       \\
       &u^{2/3}&       &       &       &       &       &       &       \\
       &       &u^{1/3}&       &       &       & 1{-}u   &       &       \\
       & 1{-}u   &       &u^{1/3}&       &       &       &       &       \\
       &       &       &       &   1   &       &       &       &       \\
u^{2/3}&       &       &       &       &u^{2/3}&       &  1{-}u  &       \\
       &       &       &       &       &       &u^{2/3}&       &       \\
  {-}1   &       &       &       &       &       &       &u^{1/3}&       \\
       &       &       &       &\phantom{u^{1/3}}&&    &       &   1
       \ema\eenon\\[-3ex]
(C1)\\*[-1ex]\ii
\benon\arraycolsep7mm\bma{ccccccccc}
\o{\quad\ \  {-}2\z_9^2{+}\z_9{-}1}&&&&&\o{{-}\z_9^8{-}\z_9^6}&&
\o{{-}\z_9^6{-}\z_9^4}&\\
&\o{{-}\z_9^4{-}\z_9^2}&&\o{\z_9^8{-}\z_9^7{+}\z_9^6}&&&&&
\o{{-}\z_9^6{-}\z_9^4}\\
&&\o{\z_9^3}&&\o{\z_9^7}&&\o{\z_9^8{-}\z_9^7{-}\z_9^2{-}1}&&\\
&\o{\z_9^6{+}\z_9^4{-}\z_9^2}&&\o{{-}\z_9^4{-}\z_9^2}&&&&&
\o{{-}\z_9^5{-}\z_9^3}\\
&&\o{{-}\z_9^8{+}\z_9^7{+}1}&&\o{{-}\z_9^7{+}\z_9^6{+}\z_9^4}&&
\o{\z_9^5{-}\z_9^4{+}\z_9^3}&&\\
\o{\quad {-}\z_9^8{+}\z_9^7{+}1}&&&&&\o{\z_9^4{-}\z_9^3{+}\z_9^2}&&
\o{\z_9^4}&\\
&&\o{{-}\z_9^2{-}1}&&\o{\z_9^2}&&\o{\z_9^3}&&\\
\o{\ \ \z_9^2{-}\z_9{+}1}&&&&&\o{\z_9^8{-}\z_9^7{+}\z_9^6{+}\z_9^4{-}1}&&
\o{\z_9^4{-}\z_9^3{+}\z_9^2}&\\
&\o{\z_9}&&\o{\z_9^2}&&&&&\o{{-}\z_9^7{+}\z_9^3{-}\z_9^2{-}1\qquad\quad}
\ema\eenon
\\
(D1/D2)\\*[-6ex]\ii
\benon\bma{ccccccccc}
   1   &       &       &       &       &       &       &       &       \\
       &\g_3^2\z_3&    &       &       &       &       &       &       \\
       &       &\g_3\z_3^2&    &       &       &       &       &       \\
       &1{-}\z_3^2&      &\g_3\z_3&      &       &       &       &       \\
       &       &{-}\si\g_3\z_3^2&   &   1   &       &   \si &       &       \\
\g_3^2\z_3^2&  &       &       &       &\g_3^2\z_3^2& &1{-}\z_3^2 &     \\
       &       &1{-}\z_3^2&      &       &       &\g_3^2 &       &       \\
{-}1     &       &       &       &       &       &       &\g_3   &       \\
       &       &       &       &       &       &       &       &   1
\ema\eenon
\rightline{$\si=0$ (D1) or 1 (D2)}\\[-2ex]
(E1/E2)\\*\ii
\benon\bma{ccccccccc}
   1   &       &       &       &       &\si\frac 13(\z_3{-}1)& &
\si\tfrac13(\z_3^{1/3}{-}\z_3^{4/3})&\\
       &\z_3^{5/3}& &1{-}\z_3^2 &&       &       &       &\si\z_3^{5/3}\\
       &       &\z_3^{1/3}&    &{-}1     &       &       &       &       \\
       &       &       &\z_3^{1/3}&    &       &       &       &{-}\si \\
       &       &       &       &   1   &       &       &       &       \\
       &       &      &       &       &\z_3^{8/3}&     &       &       \\
       &       &1{-}\z_3^2&      &\z_3^{5/3}&    &\z_3^{5/3}&    &       \\
       &       &       &       &       &1{-}\z_3^2&      &\z_3^{7/3}&    \\
       &       &       &       &       &       &       &       &   1
\ema\eenon \rightline{$\si=1$ (E1) or 0 (E2)}

\newpage

(F1/F2)\\*[-6ex]\ii
\benon\arraycolsep2.5mm\bma{ccccccccc}
   1   &\m\frac12&      &\frac12&\frac14&       &       &       &       \\
       &1      &       &       &\m\frac12&       &       &       &       \\
       &       &1      &       & \m\si &       &       &       &       \\
       &       &       &1      &\frac12&       &       &       &       \\
       &       &       &       &   1   &       &       &       &       \\
       &       &       &       &       &   1   &       &       &       \\
       &       &       &       &  \si  &       &   1   &       &       \\
       &       &       &       &       &       &       &1      &       \\
       &\m 2\nu&       &  2\nu & \nu   &       &       &       &   1
\ema\eenon \rightline{$\si=0$ (F1) or 1 (F2)}

(F3/F4)\\*[-6ex]\ii
\benon\arraycolsep2.5mm\bma{ccccccccc}
   1   &\m\frac16&      &\frac16&\frac1{18}&       &       &       &       \\
       &1      &       &       &\m\frac13&       &       &       &       \\
       &       &1      &       & \m\si  &\frac16&       &\m\frac12&      \\
       &       &       &1      &\frac13&       &       &       &       \\
       &       &       &       &   1   &       &       &       &       \\
       &       &       &       &       &   1   &       &       &       \\
       &       &       &       &  \si  &\frac12&   1   &\m\frac16&      \\
       &       &       &       &       &       &       &1      &       \\
       &\m 2\nu&       &  2\nu &\frac23\nu&    &       &       &   1
\ema\eenon \rightline{$\si=0$ (F3) or 1 (F4)}

(F5/F6)\\*[-6ex]\ii
\benon\bma{ccccccccc}
   1   &{-}1{-}\nu&       &{1{+}\nu}&{-}\nu(1{+}\nu)&   &       &  &   \\
       &1      &       &       & \nu   &       &       &       &       \\
       &       &1      &       &\m \si &{-}1{-}2\nu&       &       &   \\
       &       &       &1      &\m\nu  &       &       &       &       \\
       &       &       &       &   1   &       &       &       &       \\
       &       &       &       &       &   1   &       &       &       \\
       &       &       &       &  \si  &       &   1   &{1{+}2\nu}&    \\
       &       &       &       &       &       &       &1      &       \\
       &       &       &       &       &       &       &       &   1
\ema\eenon \rightline{$\si=0$ (F5) or 1 (F6)}

(G2)\\*[-6ex]\ii
\benon\bma{rrrrrrrrr}
   1   &{{-}1}   &   1   &1      &       &       &       &       &       \\
       &1      &{{-}1}   &       &       &       &       &       &       \\
       &       &1      &       &       &{{-}1}   &       &       &1      \\
       &       &       &1      &       &       &   1   &  1    &       \\
       &       &       &       &   1   &       &       &       &       \\
       &       &       &       &       &   1   &       &       &{{-}1}   \\
       &       &       &       &       &       &   1   &1      &       \\
       &       &       &       &       &       &       &1      &1      \\
&&&\phantom{{-}1}&\phantom{{-}1}&&\phantom{{-}1}        &\phantom{{-}1}& 1
\ema\eenon

\newpage
\begin{samepage}
\subsection*{Acknowledgments}
We thank D.~I.~Gurevich for introducing us to his results,
and R.~Coqueraux for comments on an early version of this paper.
H. E. is grateful to M.~Niedermaier and S.~Theisen for fixing the
English and some logic in this paper.

\end{samepage}
\end{document}